# Decision Making: Superiority Degree

Vladimer Zhukovin, Zurab Alimbarashvili

**Abstract**
It is introduced the concept of Superiority Degree one competitive decision over another. On the basis of this concept the mathematics – theoretic structure is developed, which is part of pairs - comparisons branch in modern decision making theory. It will be useful for practice and interesting for scientific research.

**Key words:** Decision making, Pairs comparison, Superiority Degree, Incomplete information.

## l. Introduction.

The famous American cybernetics F.George, analyzing in his book: "The Foundations of Cybernetics" the stages of cybernetics' development in future, writes: That we may expect progress in all schools, but essentially in creation of the computer programs keeping decision making principles on the analogy of the same ones in humans.

The correspondent investigation on decision making has been conductor in the Institute of Cybernetics of Georgia from the very beginning of its foundation. As this paper would has been writing for the collection devoted to $40^{th}$ jubilee of our Institute, let us turn back to history.

The year 1967, The Third All - Union Symposium on Cybernetics organized by the Institute of Cybernetics of Georgia - "Decision making in humans" was one of the three problems being discussed. The thesis's book has been published.

The year 1972, The Sixth All - Union Symposium on Cybernetics organized by the Institute of Cybernetics of Georgia - devoted clear to "Decision making in humans' all kind activities. Six volumes of thesis were published. (Russian).

The year 1979, The Soviet-American Symposium on the subject of Normative and Descriptive Models of decision making, Tbilisi, the Institute of Cybernetics of Georgia was one among the organizers. The book of accounts named similarly has been published(English and Russian).

The year 1985, Berlin, Symposium on System Analysis and Simulation. Many accounts devoted to the decision making systems (complex) opposite decision making procedures were presented. The plenary report has been suggested from the Institute of Cybernetics of Georgia.

There were other International Fo>:ums included always the decisions making sections. Every 2-3 years during the $70-80^{th}$ the Seminar on the Operation Research and System Analysis has been conducting by the laboratory of "Decision Making Theory", very popular among the decision making specialists as it was the sole in USSR wholly devoted to multicriteria decision making problems. It was conducted in Kutaisi or Batumi as usual and was named by ORSA-N.

Overlooking the numerous of our publications including foreign ones, note that the most resonance has been attended in the main to our investigations on Fuzzy Multicriteria Decision Making Problems. However we have at once interesting results in other branch of the Decision making theory. There are investigations associated with concept of the superiority degree of one decision over another. Such conception has been introduced at first in [1] in context of group (social) decisions of the modern Decision making theory which are very wide and multiform. The pair-comparison branch of this theory is one of the most important being the basic for other ones. We have presented in this paper the theoretic-mathematical construction (the ABC of theory) of the degree superiority for the pair-comparison.

## 2. Preference Relations

In scientific research on decision making problems preference relations occupy the most part as research tools. They are defined and noted first by Fishborn [3]. It is no wonder because they allow to compare one decision with other so as to choose better of them in some concrete situation of decision making. Theory of preference relations is elaborated now already [2,7]. This theory, which is necessary for this article.

Let X is set of competitive decisions (alternatives), it is finite. Then $E = X \times X$ is set of all ordered pairs of decisions. Mathematically preference relation is defined as $R \subseteq E$, i.e. it is also set and all set theory operations may be used on $R$. Any preference relation is binary relation. Inverse preference relation $R^{-1}$ corresponds to any $R$ in the following way: if pair $(x, y) \in R$, then pair $(y, x) \in R^{-1}$. Inverse preference relation $R^{-1}$ has all those characteristics which corresponding preference relation $R$ has. Any $R$ is composed by two components: identity relation $R^e = R \cap R^{-1}$, strict preference relation $R^S = R \setminus R^{-1}$ and, it is clear, that $R = R^e \cup R^s$, but $R^e \cap R' = \emptyset$. If each of them is transitive, then $R^e$ is equivalence relation and $R^s$ is strict order. One of them may be empty. But if both preferences are non-empty, then $R$ is non-strict preference relation or quasi-order when it is transitive. Any preference relation $R$ may be connected or disconnected. It is connected when all decisions pairs from $E$ are comparable by $R$, i.e. $(x\ y) \in R$ or $(y, x) \in R^{-1}$ are fulfilled simultaneously. Some set $X_\Pi(R) \subseteq X$ corresponds to preference relation $R$. It is its core, which contains maximal or effective decisions. We name it Pareto-set similarly with multicriteria decision making problems. It may be empty, but if $R$ is transitive then always $X_\Pi(R) \neq \emptyset$. Let present also following results:

1. For two preference relations $R_1$ and $R_2$ if condition $R_1^s \subseteq R_2^s$, then $X_\Pi(R_2) \subseteq X_\Pi(R_1)$.

2. Preference relation $R_1$ is coordinated with preference relation $R_2$ if conditions $R_1^s \subseteq R_2^s$ and $R_1^\ell \subseteq R_2^\ell$ fulfill simultaneously. In this case $X_\Pi(R_2) \subseteq X_\Pi(R_1)$.

We shall present the other results and notions when they will be needed.

2. **Superiority Degree (SD).** Let present two classes of functions:

$$H = \{\psi(x, y) | \psi(x, y) = -\psi(y, x)\}$$

They are skew - symmetric functions.

$$T = \{\varphi(x, y) | \varphi(x, z) + \varphi(z, y) = \varphi(x, y)\}$$

this is a condition of transitivity.

Functions $\psi(x, y)$ and $\varphi(x, y)$ are scalar and x, y, z are from X - set.

**Definition 1.** Let note any scalar function $\varphi(x, y)$ defined on the set E as **superiority degree** of one competitive decision $x \in X$ over other competitive decision $y \in X$ if it is skew-symmetric, i.e. $\varphi(x, y) \in H$.

**Basic characteristics of SD:**

a) $-\varphi^* \leq \varphi(x, y) \leq \varphi^*$ , where $\varphi^* = \max_{(x, y) \in E} \varphi(x, y)$

b) $\varphi(x, y) + \varphi(y, x) = 0$

c) $\varphi(x, x) = 0$. If $\varphi(x, y) = 0$, then x and y are identical

d) $\sum_{x \in X} \sum_{y \in X} \varphi(x, y) = 0$ (3)

$$\sum_{x \in X} \sum_{y \in X} \lambda(x)\lambda(y)\varphi(x, y) = 0$$

where $\lambda(x) \geq 0$ and $\sum_{x \in X} \lambda(x) = 1$

The coefficients $\lambda(x)$ are interpreted as "weights" or coefficients of significance of decisions.

In contrast to rations scale all results in this article, stated by. using SD, will be given in differences scale. We make a remark because all known publications on these problems use rations scale.

**Definition 2.** Let note following scalar function $F(x, y)$ defined on the set $E$ as Integral (global) Superiority Degree (**ISD**):

$$F(x, y) = \sum_{z \in X} \lambda(z) \cdot [\varphi(x, z) - \varphi(y, z)] \quad (4)$$

where $\varphi(x, y)$ is **SD**, i.e. $\varphi(x, y) \in H$.

In this case two decisions x and y are compared with one an other non-directly but by means of third decision (reference point) z. Because $\varphi(x, y)$ is difference estimation the difference in formula *4* don't contradict to common sense. **Basic characteristics of ISD:**

  *1.* $F(x, y) = -F(y, x)$, i.e. $F(x, y) \in H$ and it is superiority degree.

  2. Hence it possesses all characteristics of **SD** (formulae 3),

  3. $F(x, y) \in T$ always independently of corresponding characteristic for $\varphi(x, y)$,

  4. But if $\varphi(x, y) \in T$, then $F(x, y) = \varphi(x, y)$.

The last two characteristics are most significant. **ISD** possesses many positive aspects, which will be presented in this article. But it possesses negative aspect too. Crossing out or addition decisions in X in general influence on comparison of decisions. We take into consideration this fact and try to neutralize it by control actions.

  3. **Interconnection SD and ISD with preference relation R.** Let given any preference relation *R*. Now we can form several **SD** connected with this *R*.

This connection is based on following concept:

  **Definition 3.** Scalar function $\varphi(x, y) \in H$ is named as coordinated with preference relation *R* if the conditions:

$(x, y) \in R^s \rightarrow \varphi(x, y) > 0$ ,

$(x, y) \in R^\ell \rightarrow \varphi(x, y) = 0$ (5)

These conditions may be used for **ISD** too, i.e. for function $F(x, y)$.

On other hand, if initial information received by experts inquiry or by any other way is presented as **SD**, the corresponding preference relation may be formed always by following formula:

$$R(\ell) = \{(x, y) \in E | \varphi(x, y) > \ell\} \quad (6)$$

where constant $\ell \geq 0$. This is a binary preference relation, which we name as $\ell$-**level preference relation** and analyze it later.

  **Affirmation 1.** If some preference relation *R* is coordinated with R(0); given by formula *6,* then it is coordinated with $\varphi(x, y)$, too (definition *3* ).

Proofs of affirmations we don't present in this article. *R(0) is* always connected because $\varphi(x,y)$ is defined for all pair $(x, y) \in E$. Mean while $R$ may be disconnected. Therefore in formula 5 arrows are directed only to one side.

4. **Utility Function (UF).** Initial information for decision making problem is obtained by comparison of decisions pairs. Comparison means may be different.

But results are R or $\varphi(x,y)$. Since some superiority degree is connected with any preference relation one can assume, that initial information is presented as scalar skew symmetric function $\varphi(x,y)$, i.e. SD is presented. One of most difficult problem is the ordering of set X on the basis of the results of pairs comparisons. The ordering means to define utility function on X.

**Affirmation 2.** If $\varphi(x,y) \in H \cap T$ (formulae 1 and 2 ) then it can be represented as the difference $\varphi(x,y) = f(x) - f(y)$, where $f(x)$ is some potential function [4]. And following formula takes place:

$$f(x) = \sum_{y \in X} \lambda(y).\varphi(x,y)$$

In partial case we may assume, that $\lambda(y) = \frac{1}{n}$, where $n$ is number of competitive decisions in X . Potential $f(x)$ order set X and thus it is **utility function**. Many way are of conversion of pairs comparisons into UF. But they use intuition, experience, common sense and concrete of decision making situation. That is why we attach very importance this result (affirmation 2, formulae 7 ). This affirmation introduces in foregoing problem generality (universality), formal basis, completeness and it defines those conditions, which are necessary for problem solving. In practical work usually $\varphi(x,y) \notin T$, i.e. the transitive condition is infringed (formulae 2 ).

But then we can form (formulae 4 ) and use ISD with F(x,y), we know, that $F(x,y) \in H \cap T$ always fulfils. Hence using affirmation 2 we may write following formulae:

$$F(x,y) = q(x) - q(y) \qquad (8)$$

where $q(x) = \sum_{y \in X} \lambda(y).F(x,y)$

The last one is utility function defined on X. Thus we can order X always on the basis of data of pairs comparisons. Allowed transformations for $f(x)$ and $q(x)$ are linear: they don't infringe initial order given on X.

**Affirmation 3 .** If **SD** $\varphi(x,y) \in T$ (transitive condition), then $q(x) = f(x)$.

5. **Multicriteria utility.**

Multicriteria Decision Making Problems (MDMP) are the most wide-spread and very significant class of decision making problems in modern decision making theory. In this case each decision is estimated on the basis of several criteria. And then we have the set $\Phi$ of scalar functions, defined on E:

$$\Phi = \{\varphi_1(x,y),...,\varphi_j(x,y),...,\varphi_m(x,y)\} \qquad (9)$$

**Affirmation 4 .** If $\varphi_j(x,y) \in H$ for all $j = \overline{1,m}$ then $\varphi(x,y) \in H$ , where

$$\varphi(x,y) = \sum_{j=1}^{m} \lambda_j(y).\varphi_j(x,y).$$ If $\varphi_j(x,y) \in T$ for all $j = \overline{1,m}$ ,then. $\varphi(x,y) \in T$

Let $\varphi(x,y) \in H \cap T$, then on the basis of Affirmation 2 multicriteria utility function defined on X may be presented by following formulae:

$$L(x) = \sum_{y \in X} \lambda(y) \sum_{j=1}^{m} \lambda_j \cdot \varphi_j(x,y) \qquad (10)$$

Let us prove that it is the linear convolution known in multicriteria decision making problems.

$$L(x) = \sum_{j=1}^{m} \lambda_j \cdot K_j(x) \qquad (11)$$

where $K_j(x) = \sum_{y \in X} \lambda(y) \varphi(x,y) = f_j(x)$. This is a some effectiveness criterion of "win" type, moreover, $\lambda_j \geq 0, \sum_{j=1}^{m} \lambda = 1_j$. Hence L(x) is Pareto-effective convolution [2]. Pareto-set is formed by m effectiveness criteria $K_j(x)$.

6. $\ell$- **level preference relations.** They are connected with SD or ISD and is introduced by us (formulae 6). They are binary preference relations, i.e. $R(\ell) \in E$ for all allowed value of $\ell \geq 0$. If $\varphi(x,y) \in H$, then reverse preference relation is:

$$R^{-1}(\ell) = \{(x,y) \in E | \varphi(x,y) \leq -\ell\} \qquad (12)$$

If $\ell \neq 0$, then $R(\ell)$ is strict, disconnected and, in general, non-transitive, preference relation. If $\ell = 0$, then it is connected, non-strict and also non-transitive preference relation. Now we shall formulate transitivity conditions for $R(\ell)$.

**Affirmation 5.** If SD $\varphi(x,y) \in T$ (formulae 2), then $R(\ell)$ is transitive for all allowed $\ell \geq 0$. This means that for $\ell \neq 0$ it is strict, disconnected order, but for $\ell = 0$ it is linear quasi-order (or linear order).

On the basis of $\varphi(x,y)$ may be introduced identity relation $R^\ell = \{(x,y) \in E | \varphi(x,y) = 0\}$. If $\varphi(x,y) \in T$, then it is equivalence relation. Let introduce non-strict $\ell$-level preference relation:

$$Q(\ell) = R^\ell \cup R(\ell), \ell \neq 0 \qquad (13)$$

For $\varphi(x,y) \in T$ $Q(\ell)$ is non-strict $\ell$-level order. Let remark that for both cases disconnected.

**Affirmation 6.** If $\ell_2 > \ell_1$, then $R(\ell_2) \subseteq R(\ell_1)$ and $X_\Pi(\ell_1) > X_\Pi(\ell_2)$.

For $\varphi(x,y) \in T$ the Pareto-sets are non-empty. Hence we have formed the mathematical structure imbedding one into other non-empty Pareto-sets. This is conveniently for elaboration of dialogue procedures on computer ($\ell$ is control parameter).

$$X_\Pi(0) \subseteq X_\Pi(\ell_1) \subseteq X_\Pi(\ell_2) \subseteq X_\Pi(\ell^*) = X \qquad (14)$$

where $0 < \ell_1 < \ell_2 < \ell^*$. All $\ell > \ell^*$ haven't meaning.

**8. Similarity with fuzzy preference relations.**

**Affirmation 7.** $R(\ell_1) \cap R(\ell_2) = R(\ell)$, where $\ell = \max\{\ell_1, \ell_2\}$ and $R(\ell_1) \cup R(\ell_2) = R(\ell)$, where $\ell = \min\{\ell_1, \ell_2\}$. Thus all $\ell$-level preference relations, formed on the basis of the same SD $\varphi(x, y)$, are closed under operations of join and intersection.

Let introduce yet one class of functions:
$$S = \{\varphi(x,y) | \varphi(x,y) > \max\{\varphi(x,z), \varphi(z,y)\}\} \qquad (15)$$

**Affirmation 8.** If $\varphi(x,y) \in S$, then $R(\ell)$ is transitive for all allowed values of level $\ell$.

These facts are similar on some results from fuzzy sets theory. In future we want to determine more profound connection of superiority degree with decision making fuzzy problems.

**9. Example: group decisions.**

Initial data for group (social) decisions are described as $< X, N, R >$, where X is finite set of competitive decisions (alternatives); $N$ is number of experts in group (their indices are $v = 1 \div N$; $R$ is Vectors Preference Relation (VPR) obtained by experts inquiry and defined on set $X$; $R$ consists of $N$ components:
$$R = \{R^{(1)}, R^{(2)}, \ldots, R^{(v)}, \ldots, R^{(N)}\}, \qquad (16)$$
where $R^{(v)} \subseteq E$ is ordinary (usual) scalar binary preference relation, which map preference structure of expert with index $v$. But $E = X \times X$ is set of all ordered pairs of decisions. Group decisions itself is defined as result of some procedure over $R$, that is:
$$G = \Pi(R), \text{ and also } G \subseteq E. \qquad (17)$$

Symbol $\Pi$ haven't mathematical meaning, it's procedure notation. Let us introduce next function:
$$\delta_{ij}^{(v)} = \begin{cases} 1, & \text{if } (x_i, x_j) \in R^{(v)} \\ \frac{1}{2}, & \text{if } (x_i, x_j) \in R^{(v)} \& (x_j, x_i) \in R^{(v)} \\ 0, & \text{if } (x_j, x_i) \in R^{(v)} \end{cases} \qquad (18)$$

This function is defined on decisions pairs, i.e. $(x_i, x_j) \in E$ and other, where $x_i \in X$ and $x_j \in X$.

It is coordinated with preference relation $R^{(v)}$ in view of specific sense. We don't require transitivity of preference relations $R^{(v)}$, $v = \overline{1, N}$. For describing data, received fully from experts group, let us determine next function also defined on the decision pairs:
$$n_{ij} = \sum_{v}^{N} \delta_{ij}^{(v)} \qquad (19)$$

Now two known (traditional) group decisions can be determined in view of context of our article.

**Voting by majority:**
$$G_V^s = \{(x_i, x_j) \in E | n_{ij} = n_{ji}\} \qquad (20)$$

This is a strict preference relation, it isn't transitive in general case. Identities relation corresponds to it and is also non-transitive:
$$G_V^\ell = \{(x_i, x_j) \in E | n_{ij} = n_{ji}\} \qquad (21)$$

These two preference relation form voting by majority:
$$G_V = G_V^s \cup G_V^\ell, \qquad (22)$$

if they will be joined (together).

**K- procedure (правило Копленда):**

let us introduce following K - index for $x_i \in X$

$$\varphi_i = \sum_{j=1}^{m}(n_{ij} - n_{ji}) \qquad (23)$$

where m is number of decisions in X. Using it we can write following group decision K-procedure:

$$G_k = \{(x_i, x_j) | \varphi_i \geq \varphi_j\}. \qquad (24)$$

This is a linear order, it is transitive, K- procedure is more progressive than voting by majority, because each decision is compared with all other decisions from set X. But it is true that it is not ideal.

Let introduce following number function, defined on pairs of decisions:

$$Z_{ij} = n_{ij} - n_{ji}. \qquad (25)$$

It is obvious that it is SD, because $Z_{ij} = -Z_{ji}$. Let introduce ISD too:

$$F_{ij} = \sum_{s=1}^{m}(Z_{is} - Z_{js}) \qquad (26)$$

It is easy to prove that following condition takes place:

$$F_{ij} = V_i - V_j \qquad (27)$$

where $V_i = V(x_i)$ is number potential function, defined on the set X and given on difference scale. It is (formulae 8 and 13):

$$V_i = \sum_{s=1}^{m} Z_{is} = \varphi_i \qquad , \text{K-index} \qquad (28)$$

Let determinate now $\ell$ -level preference relation:

$$G(\ell) = \{(x_i, x_j) \in E | F_{ij} \geq \ell\} \qquad (29)$$

where $\ell \geq 0$. This is a disconnected, strict order when $\ell \neq 0$ and it is a linear order when $\ell = 0$. Level $\ell$ of preference relation $G(\ell)$ can be selected on the bases of practical considerations as in work [3] for example.

The interesting results can be formulated at once:
1. G(0) is K - procedure.
2. Condition $Z_{ij} = Z_{is} + Z_{sj}$ is sufficient condition for transitivity of Voting by Majority.
3. When this condition takes place then K-procedure and Voting by Majority are equivalent [6].

As soon as $\ell$-level group preference relation G($\ell$) is formed one can define $\ell$-level group decision - this is a core of G($\ell$), i.e. Pareto-set, which will be noted by $X_\Pi(\ell)$. Two interesting results can be proved connected with Pareto-set:

4. $X_\Pi(\ell) \neq \emptyset$ for any value of level $\ell$.

5. If two $\ell$ - level group preference relations are with levels $\ell_1$ and $\ell_2$ correspondingly and $\ell_1 > \ell_2$ then condition $X_\Pi(\ell_2) \subseteq X_\Pi(\ell_1)$ (affirmation 6).

Thus we can form the structure of imbedding one to an other non-empty Pareto-sets.

Let remark only that the foregoing problems presented in this article deal with complete initial information.

Now we present one variant of using of SD and ISD in decision making problems with incomplete initial information.

**10. Incomplete information: disconnected preference relation.** In practical work the situations with incomplete information arise very often when some part of decisions pairs remain incomparable. The reasons of this fact may be very different: subjective as well as objective one. Disconnected preference relations correspond to this situation in mathematics.

Let disconnected preference relation $NC \in E$ is given on the set X.

Let take also any decision $x \in X$. With respect to it the set X will be separated in two parts (two subsets): $X_1(x)$ and $X_2(x)$. First subset contains decisions, comparable with x and second subset contains decisions, in comparable with x. For any $x \in X$ following conditions holds:

$$X_1(x) \cup X_2(x) = X$$
$$X_1(x) \cup X_2(x) = \varnothing \qquad (30)$$
$$x \in X_1(x)$$

Superiority degree may be introduced also in this case.

**Definition 4.**

a) Let name as Upper Superiority Degree (USD) following value:

$$u(x, y) = \begin{cases} \varphi(x, y), & \text{if } y \text{ is comparable with } x \\ \varphi^*, & \text{if } y \text{ is incomparable with } x \end{cases}$$

b) Let name as Lower Superiority Degree (LSD) following value:

$$d(x, y) = \begin{cases} \varphi(x, y), & \text{if } y \text{ is comparable with } x \\ -\varphi^*, & \text{if } y \text{ is incomparable with } x \end{cases}$$

where $\varphi(x, y) \in H$ on the set $X_1(x)$, $\varphi^*$ is maximal value of SD $\varphi(x, y)$ - formulae (30, 1) are used.

We shall explain this definition. If decisions x and y are comparable, then value $\varphi(x, y)$ for this pair belongs to interval $[-\varphi^*, \varphi^*]$. For incomparable decisions pairs we take the extreme values of this interval: very successful and very unsuccessful one.

**Characteristic of USD and LSD:**

1. $u(x, y) > d(x, y)$ always.
2. $u(x, y) = -d(y, x)$ and $d(y, x) = -u(x, y)$,
3. $u(x,x) = d(x,x) = 0$,
4. $u(x, y) > 0$ and $d(x, y) < 0$ always.

Now let introduce following utility functions:

$$f_d(x) = \sum_{y \in X} \lambda(y) d(x, y)$$

$$f_u(x) = \sum_{y \in X} \lambda(y) u(x, y) \qquad (33)$$

where $\varphi(x, y) \in T$ on the set $X_1(x)$ - formulae (30,2) and affirmation 2 are used.

Thus now the interval estimation is given for each $x \in X$ on the basis of incomplete initial information (disconnected preference relation NS):

$$\varDelta(x) = [f_d(x), f_u(x)] \qquad (34)$$

We remark only that $f_u(x) \geq f_d(x)$, and equality correspond to the point estimation of x. Now we must order the set X on the basis of interval estimations. Such problem is studied in detail by us in publications [5,8] and we don't present them here: only a little information. The disconnected, strict order with corresponding non-empty Pareto-set is formed on X. When new additional information is received the previous intervals are transformed into the intervals, which have diminished lengths. We don't take into consideration the case of false information. When we receive complete information then the lengths of all intervals $\Lambda(x)$ will be equal to 0, and we shall deal with point estimations of $x \in X$ (affirmation 2).

Let correct the formulae 33 and 34:

$$f_d(x) = \sum_{y \in X_1(x)} \lambda(y) \cdot \varphi(x,y) - \sum_{y \in X_2(x)} \lambda(y) \cdot \varphi^* = \varphi(x) - \Lambda(x) \cdot \varphi^* \qquad (35)$$

where $\varphi(x)$ is some constant value corresponding received information. It is transformed only when additional information is obtained. And $\Lambda(x) \geq 0$ is characteristic of missing information. It tends to 0, when additional information is obtained and is equal to 0 under complete information.

Similarly we may write and discuss:

$$f_u(x) = \varphi(x) + \Lambda(x)\varphi^* \qquad (36)$$

Criteria for the estimation of missing information size can be introduced, which will be useful possibly for practical work, by this way. We present three variants:

1. $\Lambda_{mean} = \sum_{x \in X} \lambda(x) \cdot \Lambda(x)$,
2. $\Lambda_{max} = \max_{x \in X} \Lambda(x)$
3. $\Lambda_{sum} = \sum_{x \in Q} \lambda(x)$, where

They are equal to 0 under complete information. Let present yet some more several results. **Affirmation 9.**

a) If the transitivity condition $d(x,y) = d(x,s) + d(s,x)$, i.e. $d(x,y) \in T$, for all allowed x, y, s, then

$$d(x,y) = f_d(x) - f_u(y) \qquad (37)$$

b) If $u(x,y) \in T$, then

$$u(x,y) = f_u(x) - f_d(y) \qquad (38)$$

**Definition 5.**

a) Let note as Integral (global) Upper Superiority Degree (IUSD) following value:

$$U(x,y) = \sum_{r \in X} \lambda(r) \cdot [u(x,r) - d(y,r)] . \qquad (39)$$

b) Let note as Integral (global) Lower Superiority Degree (ILSD) following value:

$$D(x,y) = \sum_{r \in X} \lambda(r) \cdot [d(x,r) - u(y,r)] . \qquad (40)$$

**Affirmation 10.** Following results always take place.

$$U(x,y) = f_u(x) - f_d(y) ,$$
$$D(x \bullet Y) = f_d(x) - f_u(y) . \qquad (41)$$

Let present now the characteristic of IUSD and ILSD:
1. $U(x,y) > D(x,y)$.
2. $U(x,y) = -D(y,x)$ and $D(x,y) = -U(x,y)$.
3. $D(x,y) = d(x,y)$ if transitivity condition, i.e. $d(x,y) \in T$ for all x, y, s,
   $U(x,y) = u(x,y)$, if $u(x,y) \in T$ for all x, y, s.
4. $U(x,r) + U(r,y) = U(x,y) + U(r,r)$,
   $D(x,r) + D(r,y) = D(x,y) + D(r,r)$.

5. In general $U(r, r) \neq 0$ and $D(r, r) \neq 0$.
  $D(r, r) = -U(r, r)$ and $D(r, r) \leq 0$. If all decisions pairs are comparable, then $U(r, r) = D(r, r) = 0$. It is variant with complete information.

**11. Again about group decisions: incomplete information.** The group consists of $N$ experts. Let pick out one decisions pair $(x, y) \in E$. Vector Preference Relation (VPR) $R = \{R^{(1)}..., R^{(j)},..., R^{(N)}\}$ corresponds to each pair. Analysis of it for picked out pair will give following result:
  1. $a(x, y)$ experts have voted for x.
  2. $b(x, y)$ experts have voted for y.
  3. $p(x, y)$ expert couldn't compare x toy.

All these values are the functions determined on E on the basis of $R$, and their sum is $N$ for one pair. Now let form following values:
$d(x,y) = (a(x,y) - b(x,y)) - p(x, y)$,
$u(x, y) = (a(x,y) - b(x, y)) + p(x, y)$. (42)

First of them is LSD and second is USD. This fact can be proved. Using the affirmation 10 we have:

$$f_d(x) = \sum_{y \in X} \lambda(y) \cdot \varphi(x, y) + \sum_{y \in X} \lambda(y) \cdot p(x, y)$$
$$f_u(x) = \sum_{y \in X} \lambda(y) \cdot \varphi(x, y) - \sum_{y \in X} \lambda(y) \cdot p(x, y) \qquad (43)$$

where $\varphi(x, y) = a(x, y) - b(x, y)$ and $\varphi(x, y) \in H$.

Thus the interval $\Delta(x) = [f_d(x), f_u(x)]$ corresponds to any $x \in X$ and then we shall order the set X using these interval estimations [2,5]. Let remark also that $\varphi(x, y) = \frac{1}{2}[d(x, y) + u(x, y)]$.

**12. Conclusion.** Some very simple ideas and concepts permit us to develop and to present you the mathematics -theoretic structure (the basis of theory ) for superiority degree and connected with its problems. This structure is part of pairs comparison branch in modern decision making theory. It will be useful for practical work and interesting for scientific research. Many unsolved problems are in this field yet.